\documentclass{article}

\usepackage{PRIMEarxiv}

\usepackage[utf8]{inputenc} 
\usepackage[T1]{fontenc}    
\usepackage{hyperref}       
\usepackage{url}            
\usepackage{booktabs}       
\usepackage{amsmath}
\usepackage{amsthm}
\usepackage{amsfonts}       
\usepackage{nicefrac}       
\usepackage{microtype}      
\usepackage{tabularx}
\usepackage{multirow}
\usepackage[ruled,linesnumbered]{algorithm2e}
\usepackage{fancyhdr}       
\usepackage{graphicx}       
\graphicspath{{media/}}     

\newcolumntype{Y}{>{\centering\arraybackslash}X}

\usepackage{subfig}
\usepackage{listings}

\newcommand{\AS}{{\textbf{\textit{AS}}}}
\newcommand{\MS}{{\textbf{\textit{MS}}}}

\title{A Survey on Intelligent Iterative Methods for Solving Sparse Linear Algebraic Equations
}

\author{
  Haifeng Zou \\ 
  Shenzhen Research Institute of Big Data \\
  \texttt{zouhaifeng@sribd.cn} \\
  \And
  Xiaowen Xu \\
  Laboratory of Computational Physics, Institute of Applied Physics and Computational Mathematics \\
  \texttt{xwxu@iapcm.ac.cn} \\
  \And
  Chen-Song Zhang \\
  Academy of Mathematics and Systems Science, Chinese Academy of Sciences \\
  \texttt{zhangcs@lsec.cc.ac.cn} \\
}

\begin{document}
\maketitle
\begin{abstract}
Efficiently solving sparse linear algebraic equations is an important research topic of numerical simulation. Commonly used approaches include direct methods and iterative methods. Compared with the direct methods, the iterative methods have lower computational complexity and memory consumption, and are thus often used to solve large-scale sparse linear equations. However, there are numerous iterative methods, parameters and components needed to be carefully chosen, and an inappropriate combination may eventually lead to an inefficient solution process in practice. With the development of deep learning, intelligent iterative methods become popular in these years, which can intelligently make a sufficiently good combination, optimize the parameters and components in accordance with the properties of the input matrix. This survey then reviews these intelligent iterative methods. To be clearer, we shall divide our discussion into three aspects: a method aspect, a component aspect and a parameter aspect. Moreover, we summarize the existing work and propose potential research directions that may deserve a deep investigation.
\end{abstract}


\section{Introduction}
In the field of scientific computation and industrial simulation, it often requires solving large-scale sparse linear algebraic equation systems of the following form:
\begin{align}
Ax = b \ , \label{eq-axb}
\end{align}
where $A\in \mathbb{R}^{n\times n}$ is a sparse matrix, $b\in \mathbb{R}^n$ is a right hand vector, and $x \in \mathbb{R}^n$ is a vector that need to be computed. There are mainly two types of methods for solving this equation: direct methods~\cite{davis2016survey} and iterative methods~\cite{saad2003iterative}. The main advantage of direct methods lies in their robustness, making them widely applied in the current industrial softwares. However, direct methods have higher computational complexity. For example, for dense matrices and sparse matrices, the computational complexities are $\mathcal{O}(n^3)$ and $\mathcal{O}(n^{2+})$, respectively, which limits the computational scale of problems from practice. While empirical evidence has shown that direct methods are usually efficient for solving problems with scales up to millions, iterative methods are the mainstream approach for problems with larger scales. 

An iterative method starts from an initial guess then iteratively improves the solution until it reaches a specified accuracy. A typical iteration usually includes an efficient \textit{Sparse Matrix-Vector multiplication} (\textbf{\textit{SpMV}}) operation, which results in a lower iteration complexity of $\mathcal{O}(n)$. The overall computational complexity of an iterative method is often $\mathcal{O}(kn)$, where $k$ is the number of iterations. Therefore, the computational efficiency of iterative methods essentially depends on the number of iterations. For a given problem, if a suitable iterative method can be used, then $k \ll  n$ may hold and the computational complexity can approximate $\mathcal{O}(n)$, which shows a great potential in large-scale computations.

However, compared with direct methods, iterative methods may  have robustness issues, which means that iterative methods may depend heavily on matrix properties, and there is no universal iterative method that can solve all problems optimally. Hence, we may have to select a suitable iterative method from a class of candidates, or decide an optimal combination over a set of optional parameters or components, when we face a concrete problem. While this might be easy for some problems, it is very intractable in general.

The occurrence of intelligent iterative methods shed a light upon a possible solution to this issue. They build a mapping from an \textit{Matrix Space} (\MS) to an \textit{Algorithm Space} (\AS). Here, an \AS\ consists of, e.g., candidate iterative methods, their components and parameters, while an \MS\ might be, e.g., a set of matrices $A^{(i)}$ with different properties that originate from time-dependent or nonlinear PDEs, see Eq. (\ref{eq-sequence})  below,
\begin{align}
\{ A^{(i)}  x^{(i)} = b^{(i)} , i=1,\dots,n_{eq} \}. \label{eq-sequence}
\end{align}
Due to the effect of time integration and nonlinearity, the equation number $n_{eq}$ can be the tens of thousands, of millions, or even more. Moreover, the matrix properties in \MS\  may differ significantly, which makes a single iterative method impossible to solve all the equations effectively.

The mapping obtained from an intelligent iterative method can then help select an appropriate iterative method from \AS, and optimize its parameters and components for each input matrix in \MS. This improves the solution efficiency from three aspects, i.e., method, component and parameter. Ideally, an intelligent iterative method thus determines an optimal solution method within \AS\ for each matrix, which has been a long-term pursuit of researchers. In particular, with the recent prosperity of deep learning, intelligent iterative methods have received an increasing attention. Due to their advantage and popularity, this paper provides a comprehensive overview about the recent developments of intelligent iterative methods for both interested researchers and beginners. 

The contributions of this paper are listed as follows:
\begin{itemize}
	\item We provide the first comprehensive overview concerning a crossover research combining iterative methods with the technique of deep learning. In particular, we formalize the notion of intelligent iterative methods as a hybridization of iterative methods and deep learning.
	\item In accordance with different optimization objectives, the study of intelligent iterative methods is categorized into three types, i.e., method selection based on classifiers, component optimization based on neural networks, and parameter auto-tuning based on regression models. For each type, we introduce its basic principles and then summarize their related work.
\end{itemize}

This paper is organized as follows. Section 2 introduces the background and fundamental principles of intelligent iterative methods. Section 3 presents an overall architecture of intelligent iterative methods. Section 4 discusses the current state of work on intelligent iterative methods. Section 5 summarizes related work and proposes possible research directions on intelligent iterative methods. Finally, Section 6 concludes the whole paper. 

\section{Background}
In general, iterative methods start from an initial solution and then iteratively improve the solution with a specified iterative format until the solution has reached an specified precision. Algorithm \ref{alg-iter} below shows a general framework of iterative methods for equation (\ref{eq-axb}).

\begin{algorithm}
	\SetKwInOut{KIN}{INPUT}
	\SetKwInOut{KOUT}{OUTPUT}
	\KIN{initial solution $x^{(0)}$, tolerance $\varepsilon$, max iteration number $K$ }
	$r^{(0)}=b-Ax^{(0)};$ 
	
	\For{k=1,\dots,K}{
		$x^{(k)} = G^{(k)} x^{(k-1)} + d^{(k)};$ 
		
		$r^{(k)}=b-Ax^{(k)};$  
		
		\If{$||r^{(k)}|| \leq \varepsilon || r^{(0)} ||$}{
			break 
		}
	}
	\caption{The general framework of iterative methods}
	\label{alg-iter}
\end{algorithm}

In Algorithm \ref{alg-iter}, $G^{(k)}$ and $d^{(k)}$ represent the iteration matrix and vector at step $k$ in the iterative format. Different iteration formats correspond to different iterative methods. Currently, the most widely used iterative methods are Krylov subspace methods~\cite{saad2003iterative}, including CG, GMRES, BiCGSTAB, etc. The computation process can be viewed as a combination of matrix and vector operations. For sparse matrices, the computational complexity of matrix and vector operations is typically $\mathcal{O}(n)$. Let $k$ be the number of iterations, then the overall computational complexity is $\mathcal{O}(kn)$. Therefore, the computational complexity of iterative methods mainly depends on the number of iterations, i.e., the convergence rate.

The convergence rate of Krylov subspace methods is related to the matrix condition number. For large-scale ill-conditioned problems, the convergence  is often slow, and preconditioning techniques are needed to accelerate the convergence. The so-called preconditioning technique involves constructing a matrix $M \in \mathbb{R}^{n \times n}$ such that the following preconditioned system has a relatively good condition number
\begin{align}
M^{-1}Ax = M^{-1}b \ , \label{eq-pre}
\end{align}
to mproves the convergence rate of the iterative method, where $M^{-1}$ is referred to the preconditioner.

For large practical problems, preconditioning techniques have become indispensable for accelerating iterative convergence. The question is how to construct a good preconditioners. On the one hand, it is required that $M$ is a good approximation of $A$. On the other hand, it is required that $M^{-1}$ is less computationally expensive. However, these two are often contradictory to a certain extent, e.g. when $M=A$, the approximation is satisfied, but the computation of $M^{-1}$ is too expensive. To balence these, it is often necessary to construct appropriate preconditioners according to the matrix properties and application characteristics.

Given a matrix $A$, there are usually two ways to construct the preconditioner. One is based on matrix splitting or matrix factorization. For example, by splitting the matrix into $A=D-L-U$, where $D$,$L$,$U$ are the diagonal, lower triangular, and upper triangular matrices of $A$. Classical iterative methods such as Jacobi and Gauss-Seidel can essentially be viewed as preconditioners with $M_{Jac}=D$ and $M_{GS}=D-L$, respectively. Similarly, the matrix incomplete LU decomposition (ILU) method can be viewed as a preconditioner $M_{ilu}=\tilde{L} \tilde{U}$, where $A \approx \tilde{L} \tilde{U}$. Such preconditioners can be represented explicitly in matrix form.

Another approach involves approximating $A^{-1}$ based on the computation, such as Geometric Multigrid (GMG) methods~\cite{hackbusch2013multi,brandt2011multigrid}, Algebraic Multigrid (AMG) methods~\cite{ruge1987algebraic,stuben2001review,xu2017}, and Domain Decomposition Methods (DDM)~\cite{dolean2015introduction,toselli2004domain}. Preconditioners of this type, denoted as $M$, are often not explicitly defined in specific forms, but are implicitly composed through a sequence of computation. Such preconditioners are widely used in practical problems.

Table \ref{tab-iter} lists these commonly used Krylov iterative methods from PETSc website\footnote{https://petsc.org/main/overview/linear\_solve\_table/}. It includes eight prevalent iterative methods and eight different preconditioners, resulting in a total of 64 combinations. When considering the components and parameters of each iterative method, the final strategy within the \AS\ could encompass thousands of variations. Moreover, there is no single method that can solve all problems efficiently. It is necessary to select appropriate iterative methods, adjust parameters and components based on the matrix properties in order to achieve high efficiency. The extensive \AS\ presents a significant challenge when iterative methods are used to solve practical problems.

The \MS\ comprises matrices with varying properties. Matrices primarily stem from the discretization of PDEs, resulting in matrices with distinct characteristics. Even for the same PDE, adjusting parameters within the PDE can generate matrices with diverse properties. For instance, in the case of the diffusion equations, 
\begin{equation}
\begin{split}
-  \nabla \cdot (\kappa \nabla  u) &= f_1, \quad x \in \Omega \, ,\\
u &= f_2, \quad x \in \partial \Omega \, ,
\end{split} \label{eq-diff}
\end{equation}
modifying the diffusion coefficient $\kappa$ can yield isotropic or anisotropic matrices. Similarly, adjusting the Poisson's ratio in linear elasticity equations, or modifying the permittivity in Maxwell's equations, leads to matrices with varying levels of computational complexity. Theoretically, the matrix space encompasses an infinite number of matrices.

Based on the input matrix from \MS, automatically selecting and optimizing the method, component, and parameter from \AS\ constitutes the objective of intelligent iterative methods.

\begin{table}[htbp]
	\centering
	\caption{The Algorithm Space (\textbf{\textit{AS}})  }
	\label{tab-iter}
	\begin{tabular}{|c|c|>{\centering\arraybackslash}m{4cm}|>{\centering\arraybackslash}m{5cm}|}
		\hline 
		\multicolumn{2}{|c|}{Methods} & Components & Parameters \\
		\hline
		\multirow{8}{*}[-5em]{Krylov} &  CG & - & tol \\
		\cline{2-4}
		& F-CG & - & tol, search directions, truncation type\\
		\cline{2-4}
		& GMRES & Gram Schmidt & tol, restart value, tolerance for happy ending \\
		\cline{2-4}
		& F-GMRES & Gram Schmidt & tol, restart value, tolerance for happy ending\\
		\cline{2-4}
		& L-GMRES & Gram Schmidt & tol, restart value, tolerance for happy ending\\
		\cline{2-4}
		& BICG & - & tol\\
		\cline{2-4}
		& BICGSTAB & - & tol, number of search direction \\
		\cline{2-4}
		& GCR & - & tol, restart value \\
		\hline 
		\multirow{8}{*}[-6.5em]{Precond} &  $\omega-$Jacobi & - & its, weight \\
		\cline{2-4}
		& Blocked Jacobi & block solver & its, block size \\
		\cline{2-4}
		& G-S & - & its \\
		\cline{2-4}
		& SSOR & - & its, relaxation factor\\
		\cline{2-4} 
		& GMG & coarsening, prolongation, restriction, smoother & its,cycle type\\
		\cline{2-4}
		& AMG & coarsening, prolongation, restriction, smoother & its, cycle type, strong threshold\\
		\cline{2-4}
		& DDM & subdomain solver & its, overlap size, number of subdomain \\
		\cline{2-4}
		& ILU & - & levels,  mat ordering type, pivot in blocks\\
		\hline  
	\end{tabular}
\end{table}

\section{Intelligent Iterative Methods}
Figure \ref{fig-ai-solver} below shows a general framework of intelligent iterative methods. The framework consists of five parts: algorithm space, matrix space, intelligent mapping (select, optimize, tune), intelligent iterative method, and data. The details of each part are described as follows:
\begin{itemize}
	\item Algorithm Space (\AS): including a variety of iterative methods, as well as their parameters and components. 
	\item Matrix Space  (\MS): including matrices discretized from various PDEs, such as diffusion equation, linear elastic equation, Maxwell equation.
	\item Select, Optimize, Tune: Utilizing machine learning or deep learning algorithms to intelligently match \AS\ and \MS, yielding intelligent iterative methods for efficient problem-solving. Intelligent mapping is not limited to selecting optimal combinations within given methods. It also encompasses parameter auto-tuning and component optimization based on input matrices. 
	\item Intelligent Iterative methods: Assisted by machine learning and deep learning algorithms, this type of iterative method adapts its approach based on the input matrix.
	\item Data: The data is the foundation of intelligent iterative methods, comprising matrices from the \MS\ and computation results based on the methods, components and parameters from \AS, which is utilized by ML and DL to train and test the model.
\end{itemize}

\begin{figure}[htbp]
	\centering
	\includegraphics[scale=0.35]{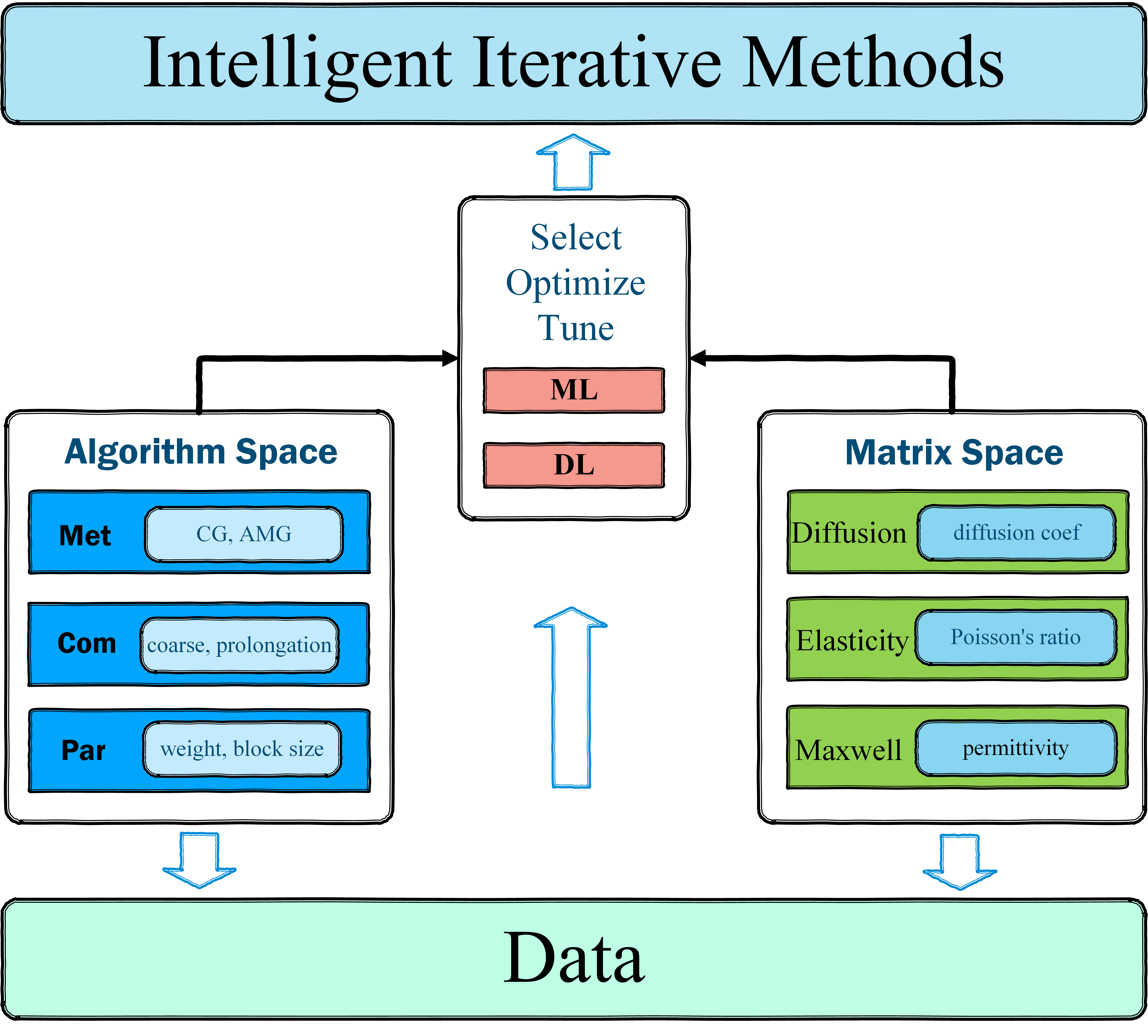}
	\caption{The framework of intelligent iterative methods, where Met stands for method, Com stands for component, Par stands for parameter}
	\label{fig-ai-solver}
\end{figure}

\section{The Development of Intelligent Iterative Methods}

Over recent decades, substantial research has been dedicated to intelligent iterative methods, exploring method, component, and parameter perspectives. Pertinent literature is outlined in Table \ref{tab-total}. Considering the varying components of different iterative methods and the diverse optimization strategies applied to distinct components, the second column ``Iterative" signifies the corresponding iterative methods encompassing these components. Meanwhile, the ``Precond" designation indicates literature employing neural networks to optimize the preconditioner $M^{-1}$. This section offers an overview of advancements and the present state of intelligent iterative methods.

\begin{table}[htbp]
	\centering
	\caption{Representative researches in three aspects}
	\label{tab-total}
	\begin{tabular}{|>{\centering\arraybackslash}m{4cm}|c|>{\centering\arraybackslash}m{7cm}|}
		\hline
		Category & Iterative  &  Publications \\
		\hline
		Methods selection & - & \cite{Eijkhout2010} \cite{Eller2012} \cite{Holloway2007}  \cite{Funk2022} \cite{Kolev2016} \cite{Motter2015} \cite{Bhowmick} \cite{Sood2018} \cite{Bhowmick2009} \cite{Sciences2016} \cite{Yamada2018} \cite{Jessup2016} \cite{tang2022graph}  \\
		\hline 
		\multirow{3}{*}{\parbox{3cm}{\centering Components \\ optimization}}&  AMG &  \cite{Katrutsa2017} \cite{Greenfeld2019} \cite{osti2021} \cite{Katrutsa2020} \cite{Chen2022} \cite{Luz2020} \cite{Huang2021} \cite{Taghibakhshi2021} \cite{Hsieh2019} \cite{Oosterlee2003} \cite{Schmitt2019}  \cite{wang4110904learning} \\
		\cline{2-3}
		& DDM & \cite{Heinlein2018} \cite{Heinlein2020} \cite{Heinlein2021} \cite{Li2020} \cite{Li2020a} \cite{sheng2022pfnn}  \cite{taghibakhshi2022learning} \cite{nguyen2022efficient} \cite{dolean2022finite} \cite{jagtap2021extended} \cite{basir2023generalized} \cite{sun2022domain} \cite{li2023deep} \\ 
		\cline{2-3}
		& Precond & \cite{Ackmann2020}  \cite{Sappl2019}  \cite{Azulai}  \cite{Luna2021} \cite{Stanaityte2020} \cite{Goetz2019} \cite{lerer2023multigrid} \cite{hausner2023neural} \cite{grementieri2022towards} \cite{gu2022deep}\\
		\hline 
		Parameters auto-tuning & - & \cite{Caldana2019} \cite{caldana2023deep}  \cite{burrows2013learning}  \cite{liu2021gptune} \\ 
		\hline 
	\end{tabular}
\end{table}

\subsection{Method Selection Based on Classifier} \label{cpt-1-3}

The method selection based on classifier involves constructing a classifier $f$ within a given set of iterative methods $C$. This approach seeks to identify the optimal method $c\in C$ based on the given matrix feature parameter $p_m$, as expressed by the equation:
\begin{equation}
c = f(p_m) \ . \label{eq-cla}
\end{equation}
The set $C$ can be made up of only preconditioners~\cite{Eller2012}, only iterative methods~\cite{Funk2022}, or the combination of preconditioner and iterative methods ~\cite{Kolev2016}.

The approach of method selection based on classifiers has long captured the attention of researchers, with its roots dating back. In 1996, Casanova and Dongarra introduced the NetSolve project~\cite{Apgar1996}, which automatically selects between iterative and direct methods based on the sparse pattern of matrices. Building upon this, Dongarra and Eijkhout initiated the SALSA project~\cite{Dongarra2003} in 2003, conducting a series of studies~\cite{Dongarra2003,EidDonEij:ipdps2003,eijkhout2003proposed,DemEtAl:ieeeproc2004,Eijkhout2010,Bhowmick}. SALSA not only automatically selects iterative methods based on matrix features but also identifies optimal sparse storage structures. Pate Motter came up with the Lighthouse Project~\cite{Motter2015,Jessup2016,Sood2018} in 2015, which compares the classification accuracy of different classifiers and different matrix features. In fact, the development of classifier-based approaches primarily focusing on matrix feature extraction and compression. 

Furthermore, method auto-selection based on classifiers is a built-in option in many commercial scientific computing and industrial software, often considered a key component of their core competitiveness. For instance, the solving command in Matlab (\lstinline|A\b|) automatically selects the appropriate method based on matrix characteristics such as sparsity, bandwidth, and data type~\cite{matlibweb}. Structural analysis software Solidworks~\cite{solidwork} can choose between Intel Direct Sparse solvers and FFEPlus Iterative method based on matrix features, while the Nastran~\cite{nastran} utilizes machine learning algorithm, Simplified Decision Trees for method auto-selection.

The essential components of method selection based on classifiers are matrix features and classification algorithms, which are illustrated below.

\subsubsection{Matrix feature}
Matrix features are pivotal factors that influence classification accuracy. Various matrix features have been developed for different problems and algorithms. In this paper, we collect features from literatures~\cite{XuLeeZhang2003,Motter2015,yue2015adaptive,xu2017algebraic,Sciences2016, eijkhout2003proposed} and categorize these features into six classes: 
\begin{itemize}
	\item Basic features: number of rows, number of nonzeros, number of nonzero diagonal elements, etc.;
	\item Structural features: symmetry, bandwidth, avergage number of nonzeros per row, number of structurally unsymmetric elements, etc.;
	\item Norm features: trace, 1 norm, infinity Norm, infinity norm of symmetric part $(A+A^T)/2$, infinity norm of anti-symmetric part $(A-A^T)/2$, etc.; 
	\item Spectral features: estimated condition number, eigenvalues, eigenvectors, singular values, etc.;
	\item Variability features: row variability 
	\begin{align*}
	\max_{i} \log_{10} \frac{\max_{j} |a_{ij}|}{\min_{j,a_{ij}\ne0}  |a_{ij}|} \ ,
	\end{align*}
	column variability
	\begin{align*}
	\max_{j} \log_{10} \frac{\max_{i} |a_{ij}|}{\min_{i,a_{ij}\ne0}  |a_{ij}|} \ ,
	\end{align*}
	average value of absolute diagonal elements, standard deviation of diagonal average, etc.; 
	\item Multiscale features: the multiscale strength of rows
	\begin{align*}
	\psi := \lfloor \log_{10} (\max_{i}v_{i})\rfloor, \quad \text{ where }
	v_{i} = \frac{\max_{j\ne i}|a_{ij}|}{\min_{j\ne i, a_{ij}\ne0}|a_{ij}|} \ ,
	\end{align*}
	the multiscale strength of columns, etc.
\end{itemize}
Noted that, the computation complexity varies across different features. For instance, basic features and structural features typically exhibit a computational complexity of $\mathcal{O}(n)$ or $\mathcal{O}(nnz)$, while the complexity of certain features can be $\mathcal{O}(1)$. In contrast, spectral features entail a higher computational complexity, often comparable to that of solving the linear equation. One of the most renowned matrix feature extraction tools is AnaMod~\cite{eijkhout2003proposed,Arnold2000} from the SALSA project, which is a parallel feature extraction program based on PETSc~\cite{balay2019petsc}. However, it is no longer maintained.

Although different features reflect matrix properties from different perspectives, their precise influence on efficiency requires further investigation. In order to reduce the number of features, it is typically necessary to select appropriate features as the input of the classifier. For instance, Bhowmick et al.~\cite{Bhowmick2009} calculated the values of each feature in different matrices, then computed the variance of each feature, and select those with variances exceeding a specific threshold. Lighthouse Project~\cite{Motter2015,Jessup2016,Sood2018} utilized RemoveUseless filter from Weka software~\cite{hall2009weka} to select features. Jae-Seung et al.~\cite{Kolev2016} utilized the Gradient Boosting Machine (GBM) algorithm~\cite{friedman2001greedy} to compute effect values of various features, selecting those with effect values surpassing a certain threshold.

\begin{figure}
	\centering
	\begin{minipage}[t]{0.5\linewidth}
		\centering
		\includegraphics[width=2.1in]{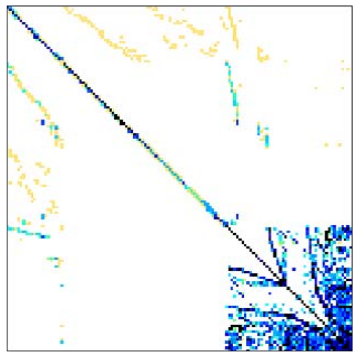}
	\end{minipage}%
	\begin{minipage}[t]{0.5\linewidth}
		\centering
		\includegraphics[width=2.1in]{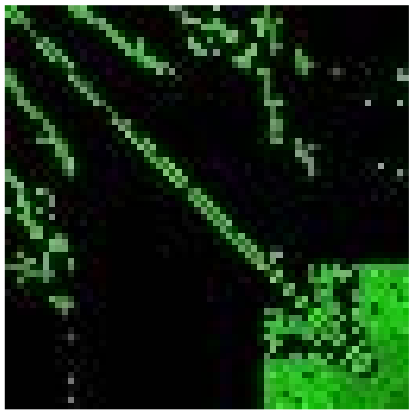}
	\end{minipage}
	\caption{Left: the sparse pattern of the matrix, right: the image transformed from the matrix (pictures are from~\cite{Yamada2018})}
	\label{fig-mp}
\end{figure}

In addition to directly extracting the aforementioned algebraic features, an alternative approach involves transforming matrices into images and subsequently employing CNN (Convolutional Neural Network) to extract image features as matrix features~\cite{Yamada2018}. As illustrated in Figure \ref{fig-mp}, sparse matrices can be mapped to a pre-defined $3\times 64 \times64$ tensor, where 3 denotes the three color channels (red, green, blue), and $64\times 64$ represents the image sizes, resulting in a colored image representation. In comparison to algebraic features, the image-based feature extraction method offers greater flexibility, allowing for arbitrary number of features.

\subsubsection{Classification Algorithm}
The currently prevalent classification algorithms encompass: Support Vector Machine (SVM)~\cite{cortes1995support}, Random Forest (RF)~\cite{breiman2001random}, Bayesian Network (BayesNet)~\cite{bielza2014discrete}, K-Nearest Neighbors (KNN)~\cite{cunningham2021k}, Alternating Decision Tree (ADT)~\cite{freund1999alternating}, C4.5 (referred to as J48 in Weka)~\cite{quinlan2014c4}, and Neural Networks.

\begin{figure}[h]
	\centering
	\subfloat[.5\linewidth][Classifying iterative methods \\ in PETSc]{
		\centering
		\includegraphics[scale=0.8]{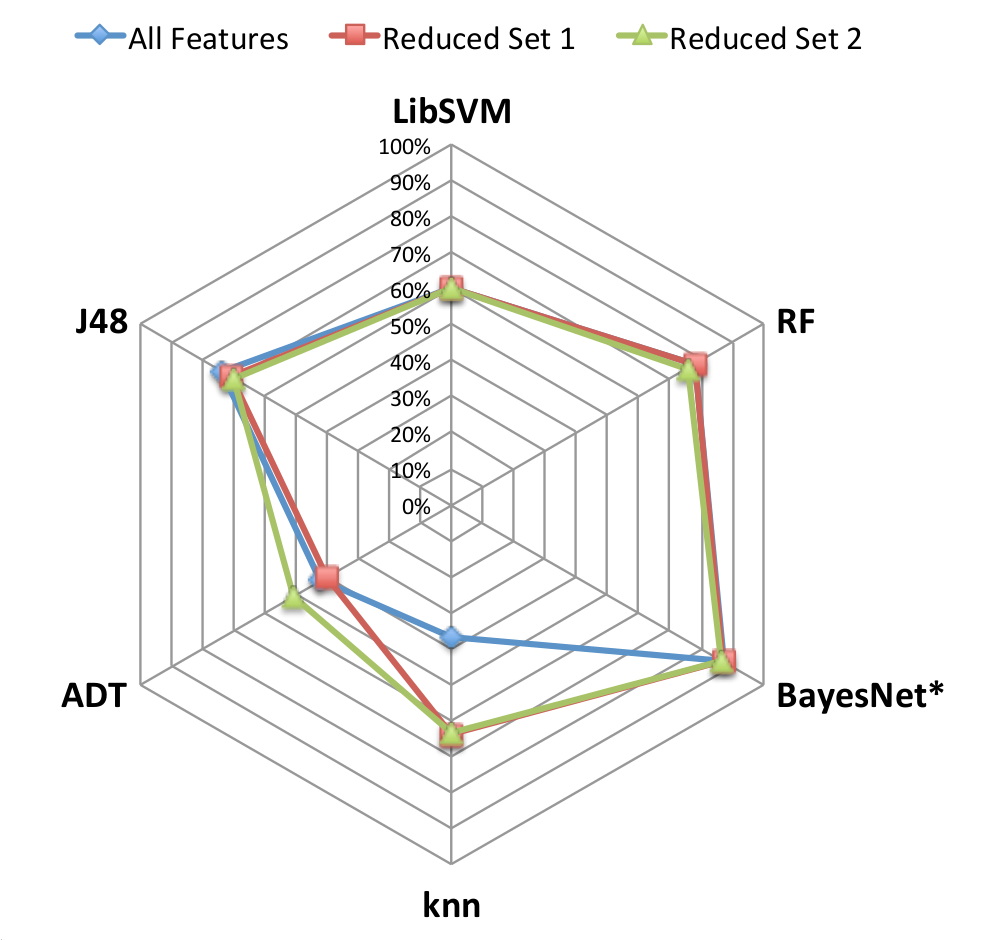}
		\label{fig-lh1}
	}
	\subfloat[.5\linewidth][Classifying iterative methods \\ in Trilinos]{
		\centering
		\includegraphics[scale=0.8]{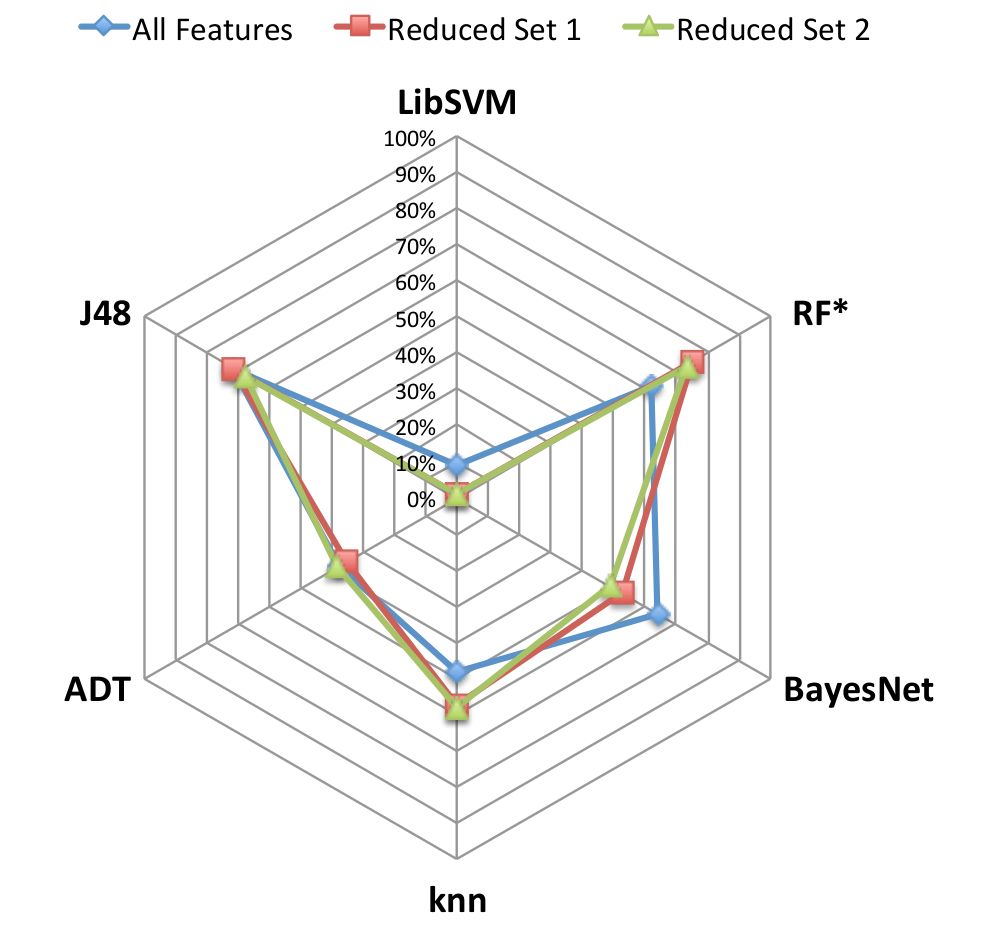}
		\label{fig-lh2}
	}
	\caption{The accuracy of the six classification algorithms based on the set $C$ composed by iterative methods from PETSc and Trilinos respectively. Three distinct colors are utilized to represent different feature sets. "All Features" encompasses all extracted matrix features, "Reduced Set 1" denotes a subset of features selected from "All Features," and "Reduced Set 2" signifies a further subset of features selected from "Reduced Set 1" (pictures are from~\cite{Motter2015}).}
	\label{fig-lh}
\end{figure}

In Lighthouse project~\cite{Motter2015,Sciences2016,Jessup2016}, iterative methods from two widely used softwares PETSc~\cite{balay2019petsc} and Trilinos~\cite{heroux2003overview} are utilized to verify the effect of the same classification algorithm based on different sets $C$. The test results, as depicted in Figure \ref{fig-lh}, indicate that the accuracy of the same classification algorithm are quite different when dealing with different set $C$. Particularly noteworthy are the SVM and BayesNet classification algorithms, where the accuracy varies by a factor of over two between the two sets. Furthermore, the optimal classification algorithms differ for the two sets, with BayesNet and RF being optimal for PETSc and Trilinos, respectively. 

Additionally, the impact of varying numbers of matrix features on classification accuracy was also investigated in the project. The results indicate that, the reduced matrix features do not yield a significant degradation on classification accuracy. On the contrary, for certain algorithms, reduced matrix features even lead to improvements in accuracy.

The following conclusions can be drawn from the data in Figure \ref{fig-lh}: Firstly, there is no universally optimal classification algorithm among all the tested algorithms. Secondly, reducing the number of matrix features not only decreases computation time but also, in some experiments, leads to improved accuracy.

\subsection{Component Optimization Based on Neural Network}

From Table \ref{tab-iter}, it can be observed that preconditioners often consist of multiple components. For instance, multigrid methods encompass smoother, grid coarsening, interpolation matrix, and restriction matrix. Similarly, domain decomposition methods involve subdomain solvers and coarse-level solvers as components. These components have diverse construction approaches, and the effectiveness of different component is closely tied to specific characteristics of the problem. In related research, component optimization is primarily based on neural networks. Different neural networks and optimization approaches are employed for distinct components. This section presents the current progress in component optimization, including multigrid (MG) and Domain Decomposition Method (DDM), etc.

\subsubsection{Mutigird Method}
The Multigrid (MG) method is a widely employed class of preconditioning acceleration methods for solving large-scale sparse linear algebraic equations in practical problems. The fundamental concept of MG involves dividing errors into high-frequency and low-frequency errors based on the eigenvalues. High-frequency errors are eliminated on fine grids, while low-frequency errors are addressed on coarse grids. This combination of coarse and fine grids results in the overall error reduction. 

MG method includes Algebraic Multigrid (AMG) and Geometric Multigrid (GMG) methods. The main distinction lies in the presence of geometric grids in GMG method, whereas AMG method don't have actual grids. Consequently, the AMG method comprises two phases: an Setup phase involving the construction of coarse grids and the corresponding coarse matrices, followed by a Solve phase wherein residual equations are solved on various grid levels. Refer to Algorithm \ref{alg-setup} and Algorithm \ref{alg-tg} for details of those two phases.

\begin{algorithm}
	Construct coarse gird according to the adjacent matrix $A$;
	
	Construct interpolation matrix $P$, and restriction matrix $R$ ($R=P^T$);
	
	Compute the coarse matrix $A_{c}=RAP$.
	
	\caption{Setup phase}
	\label{alg-setup}
\end{algorithm}

\begin{algorithm}[H]
	Pre-smoothing: smoothing $\mu_1$ times on $Ax=b$, get the approximate solution $x_f$ \\
	\eIf{deepest level }{
		Solve $Ax=b$ directly	
	} 
	{
		Restricting residuals into coarse grid: $b_c = R(b-Ax_f)$\\
		
		Solving the coarse grid equation: $A_c x_c = b_c$\\
		
		Interpolating and correcting: $x_f = x_f + Px_c$\\
	}
	Post-smoothing: smoothing $\mu_2$ times on $Ax=b$, update $x_f$
	\caption{Solve phase}
	\label{alg-tg}
\end{algorithm}
\mbox{}

The error propagation matrix $E$~\cite{baker2011multigrid} for the MG method is given by
\begin{align}
E = (I-M_2^{-1}A)^{\mu_2} (I-PA_c^{-1}RA) (I-M_1^{-1}A)^{\mu_1}, \label{eq-l1}
\end{align}
where $M_1$ and $M_2$ are the pre-smoothing and post-smoothing matrices, $\mu_1$ and $\mu_2$ are the number of pre-smoothing and post-smoothing, $A_c$ is the coarse grid matrix, and $P$ and $R$ are the interpolation and restriction matrices. A necessary and sufficient condition for the convergence of MG method is that the spectral radius of the error propagation matrix is less than 1, i.e., $\rho(E)<1$. Therefore, component optimization based on neural networks primarily focuses on relevant components that influence the spectral radius of the matrix $E$. As indicated by Eq. (\ref{eq-l1}), these components mainly include the pre-smoothing matrix $M_1$, post-smoothing matrix $M_2$, interpolation matrix $P$, and restriction matrix $R$. 

The loss is a crucial component of the neural network. For the MG method, a commonly used loss is the spectral radius of the error propagation matrix $\rho(E)$~\cite{Katrutsa2017,Greenfeld2019}:
\begin{align*}
Loss = \rho(E)  \ .
\end{align*}
However, directly computing the spectral radius $\rho(E)$ is computationally expensive. In fact, the complexity can be reduced through approximations~\cite{Katrutsa2020}. For instance, Gelfand's formula~\cite{kozyakin2009accuracy} can be employed to obtain an approximation as follows:
\begin{align*}
\rho(E) \approx \sqrt[k]{||E^k||_F},
\end{align*}
where $k$ denotes the number of iterations. Furthermore, utilizing an unbiased estimation approach~\cite{avron2011randomized}, we have
\begin{align*}
||E^k||_F^2 \approx \mathbb{E}_z(||E^kz||_2^2),
\end{align*}
where $z$ represents a vector of random variables following the Rademacher distribution, and $\mathbb{E}_z$ is the expectation. Hence, the spectral radius can be approximated as
\begin{align*}
\rho(E) \approx \left( \mathbb{E}_z(||E^kz||_2^2) \right)^{\frac{1}{2k}} 
\end{align*}
Considering that $Ez$ is equivalent to performing a single iteration to solve $Az=0$ using the MG method, then $E^kz$ corresponds to $k$ iterations, the computational complexity of the aforementioned $\rho(E)$ approximation is acceptable. For a detailed introduction, please refer to~\cite{Katrutsa2020}.

Another commonly used loss is the norm of the residuals generated during iterations~\cite{Chen2022}. Using the residuals as the loss $Loss$, the specific formula is given by
\begin{align}
Loss = || Ax_k - b||_2, \label{eq-l2}
\end{align}
where $x_k$ represents the solution at the $k$-th iteration. As shown in Algorithm \ref{alg-tg}, the computation of $x_k$ inherently involves the relevant components. Therefore, leveraging the chain rule, optimizations can be performed for the related components.

Regardless of the chosen loss, the objective is to optimize the components through neural networks. Existing works mainly fall into two categories: First, while preserving the sparse pattern of the component, optimizing the values of the components; second, directly replacing the relevant components with neural networks. Taking the interpolation matrix $P$ as an example, in the first category of methods, $P$ is equivalent to the weight coefficients in the neural network and can be represented as
\begin{align*}
Loss = f(A;P) \ ,
\end{align*}
where $A$ is the input, $Loss$ is the output, and $P$ can be denoted as the parameter to be optimized during the computation. In the second category of methods, the functionality of $P$ is replaced by the neural network $NN(w,b)$, and it can be represented as  
\begin{align*}
Loss = f(A;NN(w,b)) \ ,
\end{align*}
where NN is the Neural Network, $w$ and $b$ are the coefficients in NN. The matrix $P$ no longer exists.

The current development for these two categories are briefly outlined below.

In the first category, Katrutsa, Daulbaev, and Oseledets~\cite{Katrutsa2017} employed $\rho(E)$ as the loss function. They optimized interpolation matrix $P$, restriction matrix $R$, and the weighting coefficient $\omega$ from smoother through neural network. Specifically, each iteration step was treated as a forward propagation, and the computations within an iteration step were considered analogous to the hidden layers in the neural network. Consequently, optimizing the parameters of the neural network using backpropagation was equivalent to optimizing the matrix $P$, matrix $R$, and weighting coefficient $\omega$ in the MG method. However, the limitation of this approach is that the optimized results cannot be generalized to new, unknown matrices. Thus, for each new matrix, the matrix $P$ and $R$ must be re-trained to obtain optimized $P$ and $R$, which is impractical in reality. Greenfeld, Galun, and Kimmel~\cite{Greenfeld2019} also employed $\rho(E)$ as the loss and optimized the interpolation matrix $P$ ($R=P^T$). Differing from~\cite{Katrutsa2017}, they devised a neural network $P_{\theta}(A)$ that approximated the mapping from matrix $A$ to $P$ while preserving the sparse pattern of matrix $P$ unchanged. For an unknown matrix $A$, the trained mapping $P_{\theta}(A)$ directly yielded the optimized $P$ matrix, which enhanced the  practicality. Luz, et al~\cite{Luz2020}, Nicholas, et at.~\cite{osti2021} and Wang, et al.~\cite{wang4110904learning} presented a similar approach to~\cite{Greenfeld2019}, but with the distinction that $P_{\theta}(A)$ was constructed by the graph neural network~\cite{Battaglia2018}, improving its versatility  to non-structured grids.

In the second category, He and Xu~\cite{He2019} introduced the MgNet framework, drawing an analogy between the stencil structure generated in discrete processes and convolutional kernels in CNNs. They presented an interpretation of iterative methods from the perspective of CNNs. Building upon this, Chen and Xu~\cite{Chen2022} further proposed the Meta-MgNet algorithm. This algorithm replaces the smoothing matrix with a neural network and utilized Meta-Learning~\cite{hospedales2020meta} approach to generate the smoothing matrix based on the input matrices adaptively. Huang, Li, and Xi~\cite{Huang2021} utilized $|| x_k - x_*||_2$ as the loss, where $x_*$ is the true solution of linear equation. They used the CNN network to approximate the inverse of the smoothing matrix in AMG, corresponding to $M_1^{-1}$ in Eq. (\ref{eq-l1}). Taghibakhshi, et al~\cite{Taghibakhshi2021} employed the Dueling Double DQN algorithm from reinforcement learning~\cite{van2016deep,wang2016dueling} and the TAGCN algorithm from graph neural networks~\cite{du2017topology} to optimize the coarsening process in AMG. This optimization leads to the generation of better coarse grids, enhancing the overall computational efficiency.

Noted that as early as 2003, Oosterlee, et al~\cite{Oosterlee2003} employed a genetic optimization algorithm to select the optimal smoothing component, number of smoothing, cycle type, and parameters related to interpolation and restriction matrices in AMG. In 2019, Schmitt,  et al.~\cite{Schmitt2019} utilized genetic programming~\cite{koza1994genetic} from evolutionary computation to automatically select the smoothing component and the number of smoothing for each level of the grid in GMG. The aforementioned works can all be considered as automated optimization of components.

\subsubsection{Domain Decomposition Method}

The Domain Decomposition Method (DDM) is a preconditioner with excellent parallel scalability, developed to effectively harness the increasing computational power of supercomputers for larger-scale problems. The fundamental principle of DDM involves partitioning the computational domain into multiple subdomains and imposing appropriate boundary conditions on the subdomain interfaces. This transformation allows the problem to be decomposed into parallelizable sub-problems that can be solved concurrently on the subdomains. Depending on factors such as subdomain partitioning, internal boundary conditions, and the order of subdomain to be solved, DDM can be divided into overlapping and non-overlapping, additive and multiplicative methods, etc. Additionally, DDM can be integrated with MG method to construct multi-level DDM. Furthermore, DDM can also be used as the smoothing operator in MG method.

\begin{figure}[htbp]
	\centering
	\includegraphics[scale=1.4]{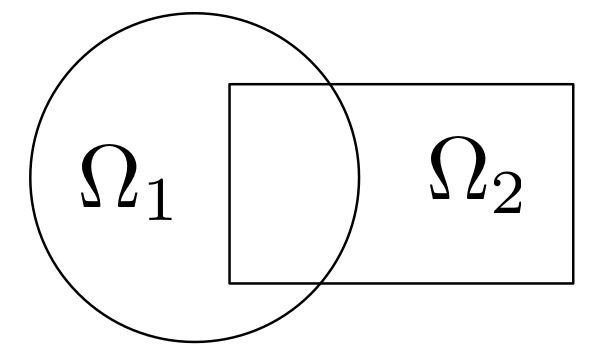}
	\caption{Two subdomains in DDM (picture is from~\cite{dolean2015introduction})}
	\label{fig-ddm}
\end{figure}

The detail of DDM is illustrated using the example of the Poisson equation:
\begin{equation}
\begin{split}
-\Delta u &= f \quad u \in \Omega ,\\
u &= g \quad  u \in \partial \Omega .
\end{split} \label{eq-poisson}
\end{equation}
The computational domain is divided into two subdomains $\Omega = \Omega_1 \cup  \Omega_2$ with an intersection $\Omega_1 \cap  \Omega_2$ as shown in Figure \ref{fig-ddm}. The boundaries are defined as $\Gamma_1 = \partial \Omega_1 \setminus \partial \Omega$ and $\Gamma_2 = \partial \Omega_2 \setminus \partial \Omega$. In DDM, the original problem is iteratively solved by solving problems in two subdomains.
\begin{equation}
\begin{cases}
-\Delta u_1^{k+1} &= f \quad \hspace{1em} u_1 \in \Omega_1 ,\\
u_1^{k+1} &= g \quad \hspace{1em} u_1 \in \partial \Omega_1 \setminus \Gamma_1 , \\
u_1^{k+1} &= u_2^k \quad  \hspace{0.5em} u_1 \in \Gamma_1 , 
\end{cases} \quad 
\begin{cases}
-\Delta u_2^{k+1} &= f \quad \hspace{1.3em}  u_2 \in \Omega_2 ,\\
u_2^{k+1} &= g \quad \hspace{1.3em}   u_2 \in \partial \Omega_2 \setminus \Gamma_2 , \\
u_2^{k+1} &= u_1^{k+1} \quad  u_2 \in \Gamma_1 .
\end{cases} \label{eq-ddm}
\end{equation}
From the above equations, it is evident that except for exchanging information on the domain boundaries $\Gamma_1$ and $\Gamma_2$, the subdomain problems can be solved independently without communication. This feature significantly enhances the parallel efficiency of DDM. Besides, Eq. (\ref{eq-ddm}) is the basic DDM; for more advanced variant of DDM, please refer to \cite{dolean2015introduction,toselli2004domain}.

The research of the combination of DDM and deep learning has been conducted by Alexander Heinlein, et al.~\cite{Heinlein2018,Heinlein2020}, and a summary of related work is provided in their review paper~\cite{Heinlein2021}. In the paper, they outline three approaches that integrate deep learning with DDM. The first approach involves utilizing deep learning to enhance the convergence rate and computational efficiency of DDM. The second approach employs neural network such as Physics-Informed Neural Networks (PINNs)~\cite{raissi2019physics} to directly solve subdomain problems, thereby bypassing conventional discretization and solving methods. The third approach leverages the principles of DDM to facilitate distributed training of neural networks, aiming to improve training efficiency. It should be noted that the third approach does not fall within the scope of DDM research, and thus, we don't introduce details of the third approach here.

In the first approach, Heinlein, et al.~\cite{Heinlein2018,Heinlein2020} optimized a variant of DDM that designed for solving problems with strong discontinuity coefficients~\cite{klawonn2017adaptive}. This method involves solving eigenvalue problems on subdomain boundaries, but the locations of boundaries are unknown. In \cite{Heinlein2018}, a fully connected neural network is employed as a classifier to determine whether a designated boundary requires solving. The input for the fully connected neural network consists of the discontinuity coefficients on both sides of the boundary. In \cite{Heinlein2020}, the CNN is employed as the classifier, where the coefficient distributions on both sides of the boundary are treated as images, and the CNN is used for image classification. Ali, et al.~\cite{taghibakhshi2022learning} utilized Graph Convolutional Neural Networks (GCNNs)~\cite{du2017topology}
and unsupervised learning to learn optimal modifications at subdomain boundaries.

The second approach involves subdomain solving methods such as Physics-Informed Neural Networks (PINNs)~\cite{raissi2019physics}, which is a novel method distinct from traditional PDE solvers. This method obviates the need for discretization of equations and instead relies on the universal approximation theorem of neural networks~\cite{scarselli1998universal} to fit solutions to PDEs. PINNs utilize randomly sampled points within the domain without the grid. The coordinates of the sampled points serve as inputs to the neural network, and the outputs are predicted solutions of the corresponding coordinates. The loss function is constructed based on the PDEs and boundary conditions. Taking Eq. (\ref{eq-poisson}) as an example, if the neural network predicts the solution $u_{pinns}$ at the coordinate $(i,j)$ within the solution domain, the loss is given by
\begin{align*}
loss_1 = (u_{pinns} - f_{ij})^2 \ .
\end{align*}
For points $(i,j)$ lying on the boundary, the loss is computed as
\begin{align*}
loss_2 = (u_{pinns} - g_{ij})^2 \ .
\end{align*}
The overall sampling loss is formulated as
\begin{align*}
loss = \sum \left[ (u_{pinns} - f_{ij})^2 + (u_{pinns} - g_{ij})^2  \right] \ .
\end{align*}
PINNs have garnered significant attention in the field of intelligent PDE solvers. As a result, a plethora of research leveraging PINNs to solve subdomain problems in DDM has emerged~\cite{nguyen2022efficient,dolean2022finite,jagtap2021extended,basir2023generalized}. Li, et al.~\cite{Li2020a} and Li et al.~\cite{Li2020} applied PINNs-like neural network to solve subdomain problems. An alternative neural network for solving PDEs is the deep Ritz~\cite{yu2018deep} method. Building upon this, Sun, et al.~\cite{sun2022domain} introduced the compensated deep Ritz method, which are utilized to solve subdomain problems. Li, et al.~\cite{li2023deep} proposed a new Multi Fourier Feature Network (MFFNet) for solving subdomain problems.

\subsubsection{Preconditioner}

The preconditioner $M^{-1}$ essentially approximates the inverse of matrix $A$, as depicted in Eq. (\ref{eq-pre}). Considering that the primary operation of the preconditioning matrix during the iteration process is matrix vector multiplication, its computation can be succinctly summarized as inputting a vector and outputting a vector. Thus, a neural network can be employed to replace the preconditioning matrix, with the caveat that the dimensions of the neural network's input and output are identical. By training the neural network using sample data, the trained network is equivalent to the optimized preconditioner. Moreover, the neural network can also be used to construct the matrix $M^{-1}$ directly or optimize the classical preconditioners.

Azulai, et al.~\cite{Azulai} and Bar, et al.~\cite{lerer2023multigrid} approximated the preconditioning matrix with the U-Net~\cite{ronneberger2015u}, which structurally resembled the MG method. The U-Net's hierarchical structure was used to replace the multigrid method's coarse grid. Sappl et al.~\cite{Sappl2019} utilized the CNN to approximate the mapping $f : \mathbb{R}^{n\times n} \rightarrow \mathbb{R}^{n\times n} , A \rightarrow M^{-1}$. To ensure the CNN's adaptability to matrices of varying sizes, the CNN architecture got rid of of fully connected layers. To prevent the generation of excessively dense matrix, four $2 \times 2$ convolutional kernels were employed within the CNN architecture. The loss was the condition number of the matrix multiplication, $cond(Af(A))$. 

Especially, for the time dependent PDE problems, there are a sequence of linear equations to be solved. Ackmann, et al.~\cite{Ackmann2020} employed the neural network to approximate the preconditioning matrix $M^{-1}$, coupled with the GCR~\cite{smolarkiewicz2000variational} method. the neural network took the solution from the previous time step to predict the preconditioner for the next time step. Luna, et al.~\cite{Luna2021} employed a CNN to fit data pairs $(b_i, x_i)$, where $b_i$ represented the $i$-th component of the right-hand side vector and $x_i$ represented the $i$-th component of the initial value. By inputting the right-hand side $b$, the CNN predicted the initial value $x$, thereby enhancing the convergence rate of GMRES during a simulation in real-time.

There exists research works that using neural network to solve the linear algebraic equation $Ax=b$ directly. Such approaches can also be treated as preconditioners to accelerate iterative convergence. Gu, et al.~\cite{gu2022deep} employed neural network to solve large-scale sparse linear equations, while Luca et al.~\cite{grementieri2022towards} utilized graph neural network to solve $Ax=b$. However, these approaches have primarily focused on standalone solutions and have not been integrated with iterative methods.

The neural network can also be used to optimize the existed preconditioners. In order to improve the efficiency of ILU method, Rita~\cite{Stanaityte2020} utilized the U-Net to optimize the sparse structure of ILU. Goetz, et al.~\cite{Goetz2019} utilized neural networks to optimize the block size within the Block Jacobi method.

\subsection{Parameter Auto-tuning Based on Regression Model}

Parameter auto-tuning based on a regression model encompass the following steps: firstly, a regression model is established to capture the relationship between matrix features, parameters to be optimized, and the metric (such as the number of iterations or wall time). Subsequently, optimizing parameters based on the  regression model. There exists a certain resemblance between regression and classification problems in terms of the rationale. As a result, part of definitions from the classification context are retained. Let $p_m$ represent matrix features, $p_a$ denote the parameters in the iterative method, $f$ denote the regression model, and $y$ denote metric associated with iterative performance. The regression model can be formulated as follows:
\begin{align}
y = f(p_m, p_a) \label{eq-reg1}.
\end{align}
Eq. (\ref{eq-reg1}) is quite similar to Eq. (\ref{eq-cla}) in classification problem. It is worth noting that Eq. (\ref{eq-reg1}) merely establish the model without conducting parameter optimization. In contrast, Eq. (\ref{eq-cla}) encompass both modeling and optimization simultaneously.

The goal of the regression model is to get the optimal parameters $p_a^*$. Given that there may exist constraints among the parameters $p_a$, the optimization problem to be solved can be formulated as follows:
\begin{alignat}{2}
\min_{p_a} \quad & f(p_m,p_a), & \notag \\
\mbox{s.t.}\quad
& f_1(p_a) = c_1,  & \label{eq-reg-c1} \\
& f_2(p_a) \geq c_2,  & \label{eq-reg-c2} 
\end{alignat}
where Eq. (\ref{eq-reg-c1}) and Eq. (\ref{eq-reg-c2}) are constraints. A commonly used approach to optimize the model is random search. After establishing the regression model based on Eq. (\ref{eq-reg1}), while keeping $p_m$ unchanged, a sequence of parameters $p_{a}^i$ satisfying the constraints is randomly generated. The outputs $y_i$ are obtained through the regression model according to the input $p_m$ and $p_{a}^i$. Finally, the index $k$ corresponding to the minimum output $y_{\min}$ is determined, yielding the optimal parameters $p_a^*=p_{a}^k$. Considering  $p_a \in \mathbb{R}^N$, the search complexity increases with larger values of $N$.

The strong connectivity threshold $\theta$ is a crucial parameter in the AMG method, whose values can influence the quantity and quality of coarse grid points, subsequently affecting convergence rate and computational efficiency. Existing solvers often resort to empirical methods, employing default values; for instance, in 2D and 3D scenarios, $\theta$ is set to 0.25 and 0.5, respectively. Matteo, et al.~\cite{Caldana2019,caldana2023deep} conducted an automatic optimization of the strong connectivity threshold $\theta$ through a regression model. The parameter $p_a$ in the regression model (Eq. (\ref{eq-reg1})) is $\theta$, while $y$ represents the approximate convergence factor of the AMG method. Matrix features $p_m$ comprise two terms: one is $-\log_2 h$ ($h$ is the grid size), and the other is the features extracted by CNN. The matrix is compressed into a $50 \times 50$ dense matrix $\hat{V}$, which can be treated as an image. Then the CNN is employed to extract image features of $\hat{V}$ as matrix features. The approach used to extract matrix features here is similar to the method employed in~\cite{Yamada2018}.

The overlap size $\delta$ can affect both the convergence rate and the communication in DDM. When $\delta$ is small, the additional computation and parallel communication are few, but the convergence rate is slow. Conversely, increasing $\delta$ accelerates convergence but also leads to higher computational and communication costs. Thus, achieving a practical balance becomes challenging. Steven, et al.~\cite{burrows2013learning} conducted an optimization of the overlap size $\delta$ in the DDM using the regression model. The inputs of the regression model are statistical features of discretized coefficients within the overlap region. The output corresponds to the number of floating-point operations, serving as a measure of algorithm efficiency.

QR decomposition can be used to solve the least squares problem in GMRES method. The block size of QR decomposition during parallel computation significantly affects overall performance. Yang, et al.~\cite{liu2021gptune} employed Gaussian processes to construct a regression model. The inputs are some common matrix features, the number of processes, and the number of threads, etc. The output is the computation time. The regression named GPTune was utilized for optimizing the block size in QR decomposition. 

\section{Summary and Future Directions}

\subsection{The Relationship Among The Three Aspects}
The researches in the three aspects mentioned above address three different types of parameters in iterative methods: method selection based on classifier focuses on discrete parameters, component optimization based on neural networks focuses on component-type parameters (such as matrices), and parameter auto-tuning based on regression models focuses on continuous parameters. These three types of parameters are the commonly used parameters in iterative methods.

From the perspective of principles, classifiers and regression models are similar. They both take a vector $x$ as input, construct a mapping function $f$, and yield a scalar output $y=f(x)$. The difference lies in the output, where $y$ is discrete in classifiers and continuous in regression models. In contrast, the approach used to optimize components varies according to the iterative methods that components belonging to. Components can be represented as matrices or scalars. This implies that, in both classifiers and regression models, scalar parameters can be optimized through component optimization methods. 

However, there exists limitation. Taking the optimization of the parameter $\theta$ in the AMG as an example~\cite{caldana2023deep}, if $\theta$ is optimized by a neural network, then an analytical expression between the loss and $\theta$ needs to be established, so that the back propagation algorithm can be utilized to optimize the parameter. The limitation is that it's nearly impossible to establish the expression between $\theta$ and the loss $\rho(E)$ or $||Ax_k -b||_2$, let alone the following optimization. In comparison, a regression model between $\theta$ and the computation time can be approximated by the neural network. However, after training the regression model, the result is a more accurate model, rather than a optimized parameter that yields the minimal computation time.

\subsection{Data Set}
Currently, research on intelligent iterative methods lack standardized datasets, resulting in variations in the matrix data used across different studies. Matrix data for component optimization typically originates from discretized model problems, such as diffusion equation (Eq. (\ref{eq-diff})). By defining different $\kappa$ values, matrices with different properties can be generated. Similarly, matrices used in regression models are also obtained through discretization of model problems. Researches in both aspects tends to optimize components and parameters within a specific iterative method for a particular problem. However, the classifiers require to select appropriate iterative methods for matrices from unknown sources. Consequently, matrices in the data set typically come from diverse problems.

We conducted a survey of the matrix data sets used in method selection problems, and the results are presented in Table \ref{tab-class}. In the table, ``Market'' is the Matrix Market data set~\cite{boisvert1997matrix}, ``SuiteSparse'' is the SuiteSparse Matrix Collection~\cite{davis2011university}, also known as the University of Florida sparse matrix collection, and ``PDEs'' denote matrices obtained by discretizing various PDEs. As evident from the table, SuiteSparse is widely utilized in research related to method selection problems.

\begin{table}[h]
	\centering
	\caption{Matrix data set used in method selection problems, ${}^*$ represents mixed data}
	\label{tab-class}
	\begin{tabular}{|c|>{\centering\arraybackslash}m{8cm}|}
		\hline
		Data Set  &  Publications \\
		\hline
		Market & \cite{Eijkhout2010} \\
		\hline 
		SuiteSparse &  \cite{Yamada2018} \cite{Motter2015} \cite{Jessup2016} \cite{Sciences2016} \cite{Bhowmick2009}  \cite{Holloway2007} \cite{Funk2022} $\cite{Kolev2016}^{*}$ \cite{tang2022graph}\\
		\hline 
		PDEs &   \cite{Eller2012} $\cite{Kolev2016}^{*}$ \cite{Bhowmick}  \\
		\hline 
	\end{tabular}
\end{table}

Currently, the SuiteSparse Matrix Collection comprises 2893 matrices\footnote{http://sparse.tamu.edu/} originating from various domains such as electromagnetics, semiconductors, fluid dynamics, structural mechanics, etc. Considering the number of matrices in SuiteSparse Matrix Collection is limited, we investigated  data sets that include matrix generation:
\begin{itemize}
	\item Structured Matrix Market (SMart)\footnote{https://smart.math.purdue.edu/}: This data set consists of 2156 matrix generation programs based on Matlab. Some of these programs are not based on PDE discretization but based on stencil structures. For instance, consider the program \lstinline|lap2d.m|, which discretizes the 2D Poisson equation using the five-point finite difference method. However, the program does not discretize the PDE but generates a sparse matrix based on the stencil structure $[-1, -1, 4, -1, -1]$. Users can only adjust the number of grid points in the $x$ and $y$ axis.
	\item SWUFE-Math Test Matrices Library\footnote{https://github.com/Hsien-Ming-Ku/SWUFE-Math}: This library includes a dozen of matrix generation programs based on Matlab.
\end{itemize}

\subsection{Matrix Feature}
In both method selection based on classifier and parameter auto-tuning based on regression model, the extraction of matrix features has always been a important point of related research. Appropriate matrix features can significantly enhance classification accuracy and reduce errors in regression models. However, there is currently no conclusive evidence in current research to indicate which matrix features are essential. Additionally, the computational cost of some features is prohibitively high.

Presently, the existing matrix feature extraction tool is AnaMod~\cite{eijkhout2003proposed,Arnold2000}, which is not longer maintained anymore. Other than conventional algebraic features, matrices can also be transformed into images, then CNN can be employed to extract features of image as matrix features~\cite{Yamada2018,Caldana2019}. This approach breaks free from the traditional features and allows for an arbitrary number of features. More importantly, it leverages end-to-end learning capabilities of neural networks to establish a direct mapping between matrices and targets, eliminating the need for manual intervention. The limitation here is that the computation domain must be rectangular, and structural grid must be employed.

A more nature approach is to use graph neural networks for matrix feature extraction~\cite{tang2022graph}. A sparse matrix can be used as the adjacent matrix of the graph. Therefore, matrices are directly input into Graph Neural Network (GNN) without additional requirements, such as computation domain and type of grid.

\subsection{Future Directions}
We suggest the following directions for intelligent iterative methods.
\begin{itemize}
	\item \textbf{Benchmark data set}: The benchmark data set serves as the foundation of intelligent iterative methods. It offers a uniform data set for training and testing for different researches, thereby furthering the development of intelligent iterative methods. While many studies select matrices from the SuiteSparse Matrix Collection to form their data set, the number and type of the selected matrices are different, making it tough to establish fair comparisons. Consequently, the benchmark data set becomes a focal point for future research.
	\item \textbf{Intelligent iterative method based on GNN}: Compared to other neural networks, GNN integrates more seamlessly with iterative methods. In classification and regression problems, GNN can be employed to extract matrix features. In the research of component optimization, GNN can be utilized to optimize preconditioner.  
\end{itemize}

\section{Conclusion}
In this survey, we propose the concept of intelligent iterative method and conduct a overview of it. The researches about intelligent iterative method are categorized into three aspects: method selection based on classifier, component optimization based on neural network, and parameter auto-tuning based on regression model. After introducing the development details of  intelligent iterative method, the relationship among three aspects, data set, matrix feature are summarized. Finally, we suggest two future directions for intelligent iterative method.

\bibliographystyle{unsrt}
\nocite{Chapelle2011,Shan2018,Tang2018,Shan2020,Wang1993,Tavakkoli2019,Brandstetter2022,Li2020b}
\bibliography{review}

\begin{thebibliography}{100}

\bibitem{davis2016survey}
Timothy~A Davis, Sivasankaran Rajamanickam, and Wissam~M Sid-Lakhdar.
\newblock A survey of direct methods for sparse linear systems.
\newblock {\em Acta Numerica}, 25:383--566, 2016.

\bibitem{saad2003iterative}
Yousef Saad.
\newblock {\em Iterative methods for sparse linear systems}.
\newblock SIAM, 2003.

\bibitem{hackbusch2013multi}
Wolfgang Hackbusch.
\newblock {\em Multi-grid methods and applications}, volume~4.
\newblock Springer Science \& Business Media, 2013.

\bibitem{brandt2011multigrid}
Achi Brandt and Oren~E Livne.
\newblock {\em Multigrid Techniques: 1984 Guide with Applications to Fluid
  Dynamics, Revised Edition}.
\newblock SIAM, 2011.

\bibitem{ruge1987algebraic}
John~W Ruge and Klaus St{\"U}ben.
\newblock {Algebraic Multigrid}.
\newblock In {\em Multigrid methods}, pages 73--130. SIAM, 1987.

\bibitem{stuben2001review}
K.~Stüben.
\newblock {A Review of Algebraic Multigrid}.
\newblock {\em Journal of Computational and Applied Mathematics},
  128(1):281--309, 2001.
\newblock Numerical Analysis 2000. Vol. VII: Partial Differential Equations.

\bibitem{xu2017}
Jinchao Xu and Ludmil Zikatanov.
\newblock {Algebraic Multigrid Methods}.
\newblock {\em Acta Numerica}, 26:591--721, 2017.

\bibitem{dolean2015introduction}
Victorita Dolean, Pierre Jolivet, and Fr{\'e}d{\'e}ric Nataf.
\newblock {\em An introduction to domain decomposition methods: algorithms,
  theory, and parallel implementation}.
\newblock SIAM, 2015.

\bibitem{toselli2004domain}
Andrea Toselli and Olof Widlund.
\newblock {\em Domain decomposition methods-algorithms and theory}, volume~34.
\newblock Springer Science \& Business Media, 2004.

\bibitem{Eijkhout2010}
Victor Eijkhout and Erika Fuentes.
\newblock {Machine Learning for Multi-stage Selection of Numerical Methods}.
\newblock {\em New Advances in Machine Learning}, pages 117--137, 2010.

\bibitem{Eller2012}
Paul~R. Eller, Jing Ru~C. Cheng, and Robert~S. Maier.
\newblock {Dynamic linear solver selection for transient simulations using
  multi-label classifiers}.
\newblock In {\em Procedia Computer Science}, volume~9, pages 1523--1532.
  Elsevier B.V., 2012.

\bibitem{Holloway2007}
America Holloway and Tzu-Yi Chen.
\newblock Neural networks for predicting the behavior of preconditioned
  iterative solvers.
\newblock In {\em International Conference on Computational Science}, pages
  302--309. Springer, 2007.

\bibitem{Funk2022}
Yannick Funk, Markus Götz, and Hartwig Anzt.
\newblock {\em Prediction of Optimal Solvers for Sparse Linear Systems Using
  Deep Learning}, pages 14--24.
\newblock 2022.

\bibitem{Kolev2016}
Jae-Seung Yeom, Jayaraman~J Thiagarajan, Abhinav Bhatele, Greg Bronevetsky, and
  Tzanio Kolev.
\newblock Data-driven performance modeling of linear solvers for sparse
  matrices.
\newblock In {\em 2016 7th International Workshop on Performance Modeling,
  Benchmarking and Simulation of High Performance Computer Systems (PMBS)},
  pages 32--42. IEEE, 2016.

\bibitem{Motter2015}
Pate Motter, Kanika Sood, Elizabeth Jessup, and Boyana Norris.
\newblock Lighthouse: an automated solver selection tool.
\newblock In {\em Proceedings of the 3rd International Workshop on Software
  Engineering for High Performance Computing in Computational Science and
  Engineering}, pages 16--24, 2015.

\bibitem{Bhowmick}
S~Bhowmick and V~Eijkhout.
\newblock {Application of machine learning to the selection of sparse linear
  solvers}.
\newblock {\em International Journal of High Performance Computing
  Applications}, pages 1--24, 2006.

\bibitem{Sood2018}
Kanika Sood, Boyana Norris, and Elizabeth Jessup.
\newblock {Comparative Performance Modeling of Parallel Preconditioned Krylov
  Methods}.
\newblock {\em Proceedings - 2017 IEEE 19th Intl Conference on High Performance
  Computing and Communications, HPCC 2017, 2017 IEEE 15th Intl Conference on
  Smart City, SmartCity 2017 and 2017 IEEE 3rd Intl Conference on Data Science
  and Systems, DSS 2017}, 2018-Janua:26--33, 2018.

\bibitem{Bhowmick2009}
Sanjukta Bhowmick, Brice Toth, and Padma Raghavan.
\newblock Towards low-cost, high-accuracy classifiers for linear solver
  selection.
\newblock In {\em International Conference on Computational Science}, pages
  463--472. Springer, 2009.

\bibitem{Sciences2016}
Kanika Sood.
\newblock {\em iterative solver selection techniques for sparse linear
  systems}.
\newblock PhD thesis, University of Oregon, 2019.

\bibitem{Yamada2018}
Kenya Yamada, Takahiro Katagiri, Hiroyuki Takizawa, Kazuo Minami, Mitsuo
  Yokokawa, Toru Nagai, and Masao Ogino.
\newblock {Preconditioner auto-tuning using deep learning for sparse iterative
  algorithms}.
\newblock {\em Proceedings - 2018 6th International Symposium on Computing and
  Networking Workshops, CANDARW 2018}, pages 257--262, 2018.

\bibitem{Jessup2016}
Elizabeth Jessup, Pate Motter, Boyana Norris, and Kanika Sood.
\newblock {Performance-based numerical solver selection in the lighthouse
  framework}.
\newblock {\em SIAM Journal on Scientific Computing}, 38(5):S750--S771, 2016.

\bibitem{tang2022graph}
Ziyuan Tang, Hong Zhang, and Jie Chen.
\newblock Graph neural networks for selection of preconditioners and krylov
  solvers.
\newblock In {\em NeurIPS 2022 Workshop: New Frontiers in Graph Learning},
  2022.

\bibitem{Katrutsa2017}
Alexandr Katrutsa, Talgat Daulbaev, and Ivan Oseledets.
\newblock {Deep Multigrid: learning prolongation and restriction matrices}.
\newblock {\em arXiv preprint arXiv:1711.03825}, nov 2017.

\bibitem{Greenfeld2019}
Daniel Greenfeld, Meirav Galun, Ron Kimmel, Irad Yavneh, and Ronen Basri.
\newblock {Learning to optimize multigrid PDE solvers}.
\newblock In {\em International Conference on Machine Learning, ICML 2019},
  volume 2019-June, pages 4305--4316, feb 2019.

\bibitem{osti2021}
Nicholas~S Moore, Eric~C Cyr, and Christopher~M Siefert.
\newblock Learning an algebriac multrigrid interpolation operator using a
  modified graphnet architecture.
\newblock Technical report, Sandia National Lab.(SNL-NM), Albuquerque, NM
  (United States), 2021.

\bibitem{Katrutsa2020}
Alexandr Katrutsa, Talgat Daulbaev, and Ivan Oseledets.
\newblock {Black-box learning of multigrid parameters}.
\newblock {\em Journal of Computational and Applied Mathematics}, 368:112524,
  2020.

\bibitem{Chen2022}
Yuyan Chen, Bin Dong, and Jinchao Xu.
\newblock {Meta-MgNet: Meta multigrid networks for solving parameterized
  partial differential equations}.
\newblock {\em Journal of Computational Physics}, 455, oct 2022.

\bibitem{Luz2020}
Ilay Luz, Meirav Galun, Haggai Maron, Ronen Basri, and Irad Yavneh.
\newblock {Learning algebraic multigrid using graph neural networks}.
\newblock In {\em 37th International Conference on Machine Learning, ICML
  2020}, volume PartF16814, pages 6445--6455, mar 2020.

\bibitem{Huang2021}
Ru~Huang, Ruipeng Li, and Yuanzhe Xi.
\newblock {Learning optimal multigrid smoothers via neural networks}.
\newblock {\em arXiv preprint arXiv:2102.12071}, 94551:1--23, 2021.

\bibitem{Taghibakhshi2021}
Ali Taghibakhshi, Scott MacLachlan, Luke Olson, and Matthew West.
\newblock Optimization-based algebraic multigrid coarsening using reinforcement
  learning.
\newblock {\em Advances in Neural Information Processing Systems},
  34:12129--12140, 2021.

\bibitem{Hsieh2019}
Jun~Ting Hsieh, Lucia Mirabella, Shengjia Zhao, Stephan Eismann, and Stefano
  Ermon.
\newblock {Learning neural PDE solvers with convergence guarantees}.
\newblock In {\em 7th International Conference on Learning Representations,
  ICLR 2019}, jun 2019.

\bibitem{Oosterlee2003}
C.~W. Oosterlee and R.~Wienands.
\newblock {A genetic search for optimal multigrid components within a Fourier
  analysis setting}.
\newblock {\em SIAM Journal on Scientific Computing}, 24(3):924--944, 2003.

\bibitem{Schmitt2019}
Jonas Schmitt, Sebastian Kuckuk, and Harald K{\"{O}}stler.
\newblock {Optimizing Geometric Multigrid Methods with Evolutionary
  Computation}.
\newblock {\em arXiv preprint arXiv:1910.02749}, oct 2019.

\bibitem{wang4110904learning}
Fan Wang, Xiang Gu, Jian Sun, and Zongben Xu.
\newblock Learning-based local weighted least squares for algebraic multigrid
  method.
\newblock {\em SSRN 4110904}, 4 2022.

\bibitem{Heinlein2018}
Alexander Heinlein, Axel Klawonn, Martin Lanser, and Janine Weber.
\newblock {Machine Learning in adaptive domain decomposition methods –
  predicting the geometric location of constraints}.
\newblock {\em SIAM Journal on Scientific Computing}, 41:A3887----A3912, 2019.

\bibitem{Heinlein2020}
Alexander Heinlein, Axel Klawonn, Martin Lanser, Janine Weber, Alexander
  Heinlein, Axel Klawonn, Martin Lanser, and Janine Weber.
\newblock {Combining Machine Learning and Adaptive Coarse Spaces - A Hybrid
  Approach for Robust FETI-DP Methods in Three Dimensions}.
\newblock {\em SIAM Journal on Scientific Computing}, 43:S816----S838, 2021.

\bibitem{Heinlein2021}
Alexander Heinlein, Axel Klawonn, Martin Lanser, and Janine Weber.
\newblock {Combining machine learning and domain decomposition methods for the
  solution of partial differential equations—A review}.
\newblock {\em GAMM Mitteilungen}, 44(1):1--28, 2021.

\bibitem{Li2020}
Ke~Li, Kejun Tang, Tianfan Wu, and Qifeng Liao.
\newblock {D3M: A Deep Domain Decomposition Method for Partial Differential
  Equations}.
\newblock {\em IEEE Access}, 8:5283--5294, 2020.

\bibitem{Li2020a}
Wuyang Li, Xueshuang Xiang, and Yingxiang Xu.
\newblock {Deep Domain Decomposition Method: Elliptic Problems}.
\newblock {\em Mathematical and Scientific Machine Learning}, 107:269--286,
  2020.

\bibitem{sheng2022pfnn}
Hailong Sheng and Chao Yang.
\newblock Pfnn-2: A domain decomposed penalty-free neural network method for
  solving partial differential equations.
\newblock {\em arXiv preprint arXiv:2205.00593}, 2022.

\bibitem{taghibakhshi2022learning}
Ali Taghibakhshi, Nicolas Nytko, Tareq~Uz Zaman, Scott MacLachlan, Luke Olson,
  and Matthew West.
\newblock Learning interface conditions in domain decomposition solvers.
\newblock {\em Advances in Neural Information Processing Systems},
  35:7222--7235, 2022.

\bibitem{nguyen2022efficient}
Long Nguyen, Maziar Raissi, and Padmanabhan Seshaiyer.
\newblock Efficient physics informed neural networks coupled with domain
  decomposition methods for solving coupled multi-physics problems.
\newblock In {\em Advances in Computational Modeling and Simulation}, pages
  41--53. Springer, 2022.

\bibitem{dolean2022finite}
Victorita Dolean, Alexander Heinlein, Siddhartha Mishra, and Ben Moseley.
\newblock Finite basis physics-informed neural networks as a schwarz domain
  decomposition method.
\newblock {\em arXiv preprint arXiv:2211.05560}, 2022.

\bibitem{jagtap2021extended}
Ameya~D Jagtap and George~E Karniadakis.
\newblock Extended physics-informed neural networks (xpinns): A generalized
  space-time domain decomposition based deep learning framework for nonlinear
  partial differential equations.
\newblock In {\em AAAI spring symposium: MLPS}, volume~10, 2021.

\bibitem{basir2023generalized}
Shamsulhaq Basir and Inanc Senocak.
\newblock A generalized schwarz-type non-overlapping domain decomposition
  method using physics-constrained neural networks.
\newblock {\em arXiv preprint arXiv:2307.12435}, 2023.

\bibitem{sun2022domain}
Qi~Sun, Xuejun Xu, and Haotian Yi.
\newblock Domain decomposition learning methods for solving elliptic problems.
\newblock {\em arXiv preprint arXiv:2207.10358}, 2022.

\bibitem{li2023deep}
Sen Li, Yingzhi Xia, Yu~Liu, and Qifeng Liao.
\newblock A deep domain decomposition method based on fourier features.
\newblock {\em Journal of Computational and Applied Mathematics}, 423:114963,
  2023.

\bibitem{Ackmann2020}
Jan Ackmann, Peter~D. D{\"{U}}ben, Tim~N. Palmer, and Piotr~K. Smolarkiewicz.
\newblock {Machine-Learned Preconditioners for Linear Solvers in Geophysical
  Fluid Flows}.
\newblock {\em arXiv preprint arXiv:2010.02866}, 2020.

\bibitem{Sappl2019}
Johannes Sappl, Laurent Seiler, Matthias Harders, and Wolfgang Rauch.
\newblock {Deep Learning of Preconditioners for Conjugate Gradient Solvers in
  Urban Water Related Problems}.
\newblock {\em arXiv preprint arXiv:1906.06925}, 2, 2019.

\bibitem{Azulai}
Yael Azulai and Eran Treister.
\newblock {Multigrid-augmented deep learning preconditioners for the Helmholtz
  equation}.
\newblock {\em arXiv preprint arXiv:2203.11025}, 2022.

\bibitem{Luna2021}
Kevin Luna, Katherine Klymko, and Johannes~P. Blaschke.
\newblock {Accelerating GMRES with Deep Learning in Real-Time}.
\newblock {\em arXiv preprint arXiv:2103.10975}, 2021.

\bibitem{Stanaityte2020}
Rita Stanaityte.
\newblock {\em ILU and Machine Learning Based Preconditioning For The
  Discretized Incompressible Navier-Stokes Equations}.
\newblock PhD thesis, University of Houston, 2020.

\bibitem{Goetz2019}
Markus Goetz and Hartwig Anzt.
\newblock {Machine learning-aided numerical linear Algebra: Convolutional
  neural networks for the efficient preconditioner generation}.
\newblock {\em Proceedings of ScalA 2018: 9th Workshop on Latest Advances in
  Scalable Algorithms for Large-Scale Systems, Held in conjunction with SC
  2018: The International Conference for High Performance Computing,
  Networking, Storage and Analysis}, pages 49--56, 2019.

\bibitem{lerer2023multigrid}
Bar Lerer, Ido Ben-Yair, and Eran Treister.
\newblock Multigrid-augmented deep learning for the helmholtz equation: Better
  scalability with compact implicit layers.
\newblock {\em arXiv preprint arXiv:2306.17486}, 2023.

\bibitem{hausner2023neural}
Paul H{\"a}usner, Ozan {\"O}ktem, and Jens Sj{\"o}lund.
\newblock Neural incomplete factorization: learning preconditioners for the
  conjugate gradient method.
\newblock {\em arXiv preprint arXiv:2305.16368}, 2023.

\bibitem{grementieri2022towards}
Luca Grementieri and Paolo Galeone.
\newblock Towards neural sparse linear solvers.
\newblock {\em arXiv preprint arXiv:2203.06944}, 2022.

\bibitem{gu2022deep}
Yiqi Gu and Michael~K Ng.
\newblock Deep neural networks for solving extremely large linear systems.
\newblock {\em arXiv preprint arXiv:2204.00313}, 2022.

\bibitem{Caldana2019}
Matteo Caldana.
\newblock {Improving Efficiency of Algebraic Multigrid Methods through
  Artificial Neural Networks}.
\newblock 2019.

\bibitem{caldana2023deep}
Matteo Caldana, Paola~F Antonietti, and Luca Dede.
\newblock A deep learning algorithm to accelerate algebraic multigrid methods
  in finite element solvers of 3d elliptic pdes.
\newblock {\em arXiv preprint arXiv:2304.10832}, 2023.

\bibitem{burrows2013learning}
Steven Burrows, J{\"o}rg Frochte, Michael V{\"o}lske, Ana
  Bel{\'e}n~Mart{\'\i}nez Torres, and Benno Stein.
\newblock Learning overlap optimization for domain decomposition methods.
\newblock In {\em Pacific-Asia Conference on Knowledge Discovery and Data
  Mining}, pages 438--449. Springer, 2013.

\bibitem{liu2021gptune}
Yang Liu, Wissam~M Sid-Lakhdar, Osni Marques, Xinran Zhu, Chang Meng, James~W
  Demmel, and Xiaoye~S Li.
\newblock Gptune: Multitask learning for autotuning exascale applications.
\newblock In {\em Proceedings of the 26th ACM SIGPLAN Symposium on Principles
  and Practice of Parallel Programming}, pages 234--246, 2021.

\bibitem{Apgar1996}
Henri Casanova and Jack Dongarra.
\newblock {NetSolve: A Network Server for Solving Computational Science
  Problems}.
\newblock {\em SuperComputing 96: Proceedings of the 1996 ACM/IEEE Conference
  on Supercomputing}, 66(December):37--39, 1996.

\bibitem{Dongarra2003}
Jack Dongarra and Victor Eijkhout.
\newblock {Self-adapting numerical software for next generation applications}.
\newblock {\em International Journal of High Performance Computing
  Applications}, 17(2):125--131, 2003.

\bibitem{EidDonEij:ipdps2003}
Thomas Eidson, Jack Dongarra, and Victor Eijkhout.
\newblock Applying aspect-orient programming concepts to a component-based
  programming model.
\newblock In {\em Proceedings of the 17th International Parallel and
  Distributed Processing Symposium (IPDPS) April 22--26, 2003, Nice, France},
  2003.

\bibitem{eijkhout2003proposed}
Victor Eijkhout and Erika Fuentes.
\newblock A proposed standard for numerical metadata.
\newblock Technical report, Citeseer, 2003.

\bibitem{DemEtAl:ieeeproc2004}
Jim Demmel, Jack Dongarra, Victor Eijkhout, Erika Fuentes, Antoine Petitet,
  Rich Vuduc, {R. Clint} Whaley, and Katherine Yelick.
\newblock Self-adapting linear algebra algorithms and software.
\newblock {\em Proceedings of the IEEE}, 93(2):293--312, 2005.

\bibitem{matlibweb}
Matlab.
\newblock mldivide, \.
\newblock \url{https://www.mathworks.com/help/matlab/ref/mldivide.html}.

\bibitem{solidwork}
Solidworks 2022.
\newblock Automatic solver selection.
\newblock
  \url{https://help.solidworks.com/2022/english/whatsnew/c_wn2022_simulation_solvers.htm}.

\bibitem{nastran}
Nastran.
\newblock Intelligent solver switch.
\newblock
  \url{https://help.mscsoftware.com/bundle/MSC_Nastran_2021.3/page/Nastran_Combined_Book/hpc_guide/chap8/chap8.xhtml}.

\bibitem{XuLeeZhang2003}
Shu{T}ing Xu, Eun-Joo Lee, and Jun Zhang.
\newblock An interim analysis report on preconditioners and matrices.
\newblock Technical Report 388-03, University of Kentucky, Lexington;
  Department of Computer Science, 2003.

\bibitem{yue2015adaptive}
Xiaoqiang Yue, Shi Shu, Xiaowen Xu, and Zhiyang Zhou.
\newblock An adaptive combined preconditioner with applications in radiation
  diffusion equations.
\newblock {\em Communications in Computational Physics}, 18(5):1313--1335,
  2015.

\bibitem{xu2017algebraic}
Xiaowen Xu and Zeyao Mo.
\newblock Algebraic interface-based coarsening amg preconditioner for
  multi-scale sparse matrices with applications to radiation hydrodynamics
  computation.
\newblock {\em Numerical Linear Algebra with Applications}, 24(2):e2078, 2017.

\bibitem{Arnold2000}
Dorian~C Arnold, Susan Blackford, Jack Dongarra, Victor Eijkhout, and Tinghua
  Xu.
\newblock {Seamless access to adaptive solver algorithms}.
\newblock {\em SGI Users' Conference}, pages 23--30, 2000.

\bibitem{balay2019petsc}
Satish Balay, Shrirang Abhyankar, Mark~F. Adams, Steven Benson, Jed Brown,
  Peter Brune, Kris Buschelman, Emil Constantinescu, Lisandro Dalcin, Alp
  Dener, Victor Eijkhout, Jacob Faibussowitsch, William~D. Gropp, V\'{a}clav
  Hapla, Tobin Isaac, Pierre Jolivet, Dmitry Karpeev, Dinesh Kaushik,
  Matthew~G. Knepley, Fande Kong, Scott Kruger, Dave~A. May, Lois~Curfman
  McInnes, Richard~Tran Mills, Lawrence Mitchell, Todd Munson, Jose~E. Roman,
  Karl Rupp, Patrick Sanan, Jason Sarich, Barry~F. Smith, Stefano Zampini, Hong
  Zhang, Hong Zhang, and Junchao Zhang.
\newblock {PETSc/TAO} users manual.
\newblock Technical Report ANL-21/39 - Revision 3.19, Argonne National
  Laboratory, 2023.

\bibitem{hall2009weka}
Mark Hall, Eibe Frank, Geoffrey Holmes, Bernhard Pfahringer, Peter Reutemann,
  and Ian~H Witten.
\newblock The weka data mining software: an update.
\newblock {\em ACM SIGKDD explorations newsletter}, 11(1):10--18, 2009.

\bibitem{friedman2001greedy}
Jerome~H Friedman.
\newblock {Greedy function approximation: a gradient boosting machine}.
\newblock {\em Annals of statistics}, pages 1189--1232, 2001.

\bibitem{cortes1995support}
Corinna Cortes and Vladimir Vapnik.
\newblock Support-vector networks.
\newblock {\em Machine learning}, 20(3):273--297, 1995.

\bibitem{breiman2001random}
Leo Breiman.
\newblock Random forests.
\newblock {\em Machine learning}, 45(1):5--32, 2001.

\bibitem{bielza2014discrete}
Concha Bielza and Pedro Larranaga.
\newblock Discrete bayesian network classifiers: A survey.
\newblock {\em ACM Computing Surveys (CSUR)}, 47(1):1--43, 2014.

\bibitem{cunningham2021k}
Padraig Cunningham and Sarah~Jane Delany.
\newblock k-nearest neighbour classifiers-a tutorial.
\newblock {\em ACM Computing Surveys (CSUR)}, 54(6):1--25, 2021.

\bibitem{freund1999alternating}
Yoav Freund and Llew Mason.
\newblock The alternating decision tree learning algorithm.
\newblock In {\em icml}, volume~99, pages 124--133. Citeseer, 1999.

\bibitem{quinlan2014c4}
J~Ross Quinlan.
\newblock {\em C4. 5: programs for machine learning}.
\newblock Elsevier, 2014.

\bibitem{heroux2003overview}
Michael Heroux, Roscoe Bartlett, Vicki Howle~Robert Hoekstra, Jonathan Hu,
  Tamara Kolda, Richard Lehoucq, Kevin Long, Roger Pawlowski, Eric Phipps,
  Andrew Salinger, et~al.
\newblock An overview of trilinos.
\newblock Technical report, Citeseer, 2003.

\bibitem{baker2011multigrid}
Allison~H Baker, Robert~D Falgout, Tzanio~V Kolev, and Ulrike~Meier Yang.
\newblock Multigrid smoothers for ultraparallel computing.
\newblock {\em SIAM Journal on Scientific Computing}, 33(5):2864--2887, 2011.

\bibitem{kozyakin2009accuracy}
Victor Kozyakin.
\newblock On accuracy of approximation of the spectral radius by the gelfand
  formula.
\newblock {\em Linear Algebra and its Applications}, 431(11):2134--2141, 2009.

\bibitem{avron2011randomized}
Haim Avron and Sivan Toledo.
\newblock Randomized algorithms for estimating the trace of an implicit
  symmetric positive semi-definite matrix.
\newblock {\em Journal of the ACM (JACM)}, 58(2):1--34, 2011.

\bibitem{Battaglia2018}
Peter~W. Battaglia, Jessica~B. Hamrick, Victor Bapst, Alvaro Sanchez-Gonzalez,
  Vinicius Zambaldi, Mateusz Malinowski, Andrea Tacchetti, David Raposo, Adam
  Santoro, Ryan Faulkner, Caglar Gulcehre, Francis Song, Andrew Ballard, Justin
  Gilmer, George Dahl, Ashish Vaswani, Kelsey Allen, Charles Nash, Victoria
  Langston, Chris Dyer, Nicolas Heess, Daan Wierstra, Pushmeet Kohli, Matt
  Botvinick, Oriol Vinyals, Yujia Li, and Razvan Pascanu.
\newblock {Relational inductive biases, deep learning, and graph networks}.
\newblock {\em arXiv preprint arXiv:1806.01261}, pages 1--40, 2018.

\bibitem{He2019}
Juncai He and Jinchao Xu.
\newblock {MgNet: A unified framework of multigrid and convolutional neural
  network}.
\newblock {\em Science China Mathematics}, 62(7):1331--1354, 2019.

\bibitem{hospedales2020meta}
Timothy Hospedales, Antreas Antoniou, Paul Micaelli, and Amos Storkey.
\newblock Meta-learning in neural networks: A survey.
\newblock {\em arXiv preprint arXiv:2004.05439}, 2020.

\bibitem{van2016deep}
Hado~van Hasselt, Arthur Guez, and David Silver.
\newblock Deep reinforcement learning with double q-learning.
\newblock In {\em Proceedings of the Thirtieth AAAI Conference on Artificial
  Intelligence}, AAAI'16, page 2094–2100. AAAI Press, 2016.

\bibitem{wang2016dueling}
Ziyu Wang, Tom Schaul, Matteo Hessel, Hado Hasselt, Marc Lanctot, and Nando
  Freitas.
\newblock Dueling network architectures for deep reinforcement learning.
\newblock In {\em International conference on machine learning}, pages
  1995--2003. PMLR, 2016.

\bibitem{du2017topology}
Jian Du, Shanghang Zhang, Guanhang Wu, Jos{\'e}~MF Moura, and Soummya Kar.
\newblock Topology adaptive graph convolutional networks.
\newblock {\em arXiv preprint arXiv:1710.10370}, 2017.

\bibitem{koza1994genetic}
John~R Koza.
\newblock Genetic programming as a means for programming computers by natural
  selection.
\newblock {\em Statistics and computing}, 4(2):87--112, 1994.

\bibitem{raissi2019physics}
Maziar Raissi, Paris Perdikaris, and George~E Karniadakis.
\newblock Physics-informed neural networks: A deep learning framework for
  solving forward and inverse problems involving nonlinear partial differential
  equations.
\newblock {\em Journal of Computational Physics}, 378:686--707, 2019.

\bibitem{klawonn2017adaptive}
Axel Klawonn, Martin K{\"U}hn, and Oliver Rheinbach.
\newblock Adaptive coarse spaces for feti-dp in three dimensions with
  applications to heterogeneous diffusion problems.
\newblock In {\em Domain decomposition methods in science and engineering
  XXIII}, pages 187--196. Springer, 2017.

\bibitem{scarselli1998universal}
Franco Scarselli and Ah~Chung Tsoi.
\newblock Universal approximation using feedforward neural networks: A survey
  of some existing methods, and some new results.
\newblock {\em Neural Networks}, 11(1):15--37, 1998.

\bibitem{yu2018deep}
Bing Yu et~al.
\newblock The deep ritz method: a deep learning-based numerical algorithm for
  solving variational problems.
\newblock {\em Communications in Mathematics and Statistics}, 6(1):1--12, 2018.

\bibitem{ronneberger2015u}
Olaf Ronneberger, Philipp Fischer, and Thomas Brox.
\newblock U-net: Convolutional networks for biomedical image segmentation.
\newblock In {\em International Conference on Medical image computing and
  computer-assisted intervention}, pages 234--241. Springer, 2015.

\bibitem{smolarkiewicz2000variational}
PK~Smolarkiewicz and LG~Margolin.
\newblock Variational methods for elliptic problems in fluid models.
\newblock Technical report, Los Alamos National Lab.(LANL), Los Alamos, NM
  (United States), 2000.

\bibitem{boisvert1997matrix}
Ronald~F Boisvert, Roldan Pozo, Karin Remington, Richard~F Barrett, and Jack~J
  Dongarra.
\newblock Matrix market: a web resource for test matrix collections.
\newblock In {\em Quality of Numerical Software}, pages 125--137. Springer,
  1997.

\bibitem{davis2011university}
Timothy~A Davis and Yifan Hu.
\newblock The university of florida sparse matrix collection.
\newblock {\em ACM Transactions on Mathematical Software (TOMS)}, 38(1):1--25,
  2011.

\bibitem{Chapelle2011}
Olivier Chapelle and Dumitru Erhan.
\newblock {Improved Preconditioner for Hessian Free Optimization}.
\newblock {\em Neural Information Processing Systems (NIPS) Workshop 2011},
  pages 1--8, 2011.

\bibitem{Shan2018}
Tao Shan, Xunwang Dang, Maokun Li, Fan Yang, Shenheng Xu, and Ji~Wu.
\newblock {Study on a 3D Possion's Equation Slover Based on Deep Learning
  Technique}.
\newblock {\em 2018 IEEE International Conference on Computational
  Electromagnetics, ICCEM 2018}, pages 3--5, 2018.

\bibitem{Tang2018}
Wei Tang, Tao Shan, Xunwang Dang, Maokun Li, Fan Yang, Shenheng Xu, and Ji~Wu.
\newblock {Study on a Poisson's equation solver based on deep learning
  technique}.
\newblock {\em 2017 IEEE Electrical Design of Advanced Packaging and Systems
  Symposium, EDAPS 2017}, 2018-Janua:1--3, 2018.

\bibitem{Shan2020}
Tao Shan, Wei Tang, Xunwang Dang, Maokun Li, Fan Yang, Shenheng Xu, and Ji~Wu.
\newblock {Study on a Fast Solver for Poisson's Equation Based on Deep Learning
  Technique}.
\newblock {\em IEEE Transactions on Antennas and Propagation},
  68(9):6725--6733, 2020.

\bibitem{Wang1993}
J.~Wang.
\newblock {Recurrent neural networks for solving linear matrix equations}.
\newblock {\em Computers and Mathematics with Applications}, 26(9):23--34,
  1993.

\bibitem{Tavakkoli2019}
Vahid Tavakkoli, Jean~Chamberlain Chedjou, and Kyandoghere Kyamakya.
\newblock {A novel recurrent neural network-based ultra-fast, robust, and
  scalable solver for inverting a “time-varying matrix”}.
\newblock {\em Sensors (Switzerland)}, 19(18):1--20, 2019.

\bibitem{Brandstetter2022}
Johannes Brandstetter, Daniel Worrall, and Max Welling.
\newblock {Message Passing Neural PDE Solvers}.
\newblock {\em arXiv preprint arXiv:2202.03376}, pages 1--27, 2022.

\bibitem{Li2020b}
Anima Anandkumar, Kamyar Azizzadenesheli, Kaushik Bhattacharya, Nikola
  Kovachki, Zongyi Li, Burigede Liu, and Andrew Stuart.
\newblock Neural operator: Graph kernel network for partial differential
  equations.
\newblock In {\em ICLR 2020 Workshop on Integration of Deep Neural Models and
  Differential Equations}, 2020.

\end{thebibliography}

\end{document}